\documentclass  [12pt,a4paper] {article}
\usepackage{amsfonts}
\usepackage{amsmath}
\usepackage{mathrsfs}
\usepackage{graphicx}
\usepackage{amsmath,amsthm,amssymb,amscd}

\newtheorem{Thm}{Theorem}[section]
\newtheorem{Lem}[Thm]{Lemma}

\theoremstyle{remark}

\theoremstyle{claim}

\DeclareMathOperator{\Diff}{Diff}

\begin{document}

\begin{center}
{\Large \bf  Diffeomorphisms with various $C^1-$(generic-)stable properties}\\
\smallskip
\end{center}

\begin{center}

\begin{center}
Xueting Tian $^{\dagger}$
\end{center}
%\smallskip
\begin{center} Academy of Mathematics and Systems Science, Chinese Academy of Sciences, Beijing 100190, China\\
\medskip
$\&$ Departamento de Matem$\acute{\text{a}}$tica, Universidade Federal de Alagoas, Macei$\acute{\text{o}}$ 57072-090, Brazil
\end{center}
%\smallskip
\begin{center}
E-mail: tianxt@amss.ac.cn $\&$ txt@pku.edu.cn
\end{center}
\medskip

Wenxiang Sun $^*$
\end{center}
%\smallskip
\begin{center}
LMAM, School of Mathematical Sciences, Peking University, Beijing 100871, China\\
\end{center}
%\smallskip
\begin{center}
E-mail: sunwx@math.pku.edu.cn
\end{center}

\footnotetext
{$^{\dagger}$ Tian is the corresponding  author and supported by CAPES.}
\footnotetext {$^*$ Sun is supported by National Natural Science
Foundation ( \# 10671006, \# 10831003) and National Basic Research
Program of China(973 Program)(\# 2006CB805903) }
 \footnotetext{AMS Review: 37A25, 37B20, 37C50, 37D20, 37D30;  }

\def\abstractname{\textbf{Abstract}}

\begin{abstract}\addcontentsline{toc}{section}{\bf{English Abstract}}
Let $M$ be a smooth compact manifold and $\Lambda$ be a compact
invariant set. In this paper we prove that for every robustly
transitive set $\Lambda$, $f|_\Lambda$ satisfies a $C^1-$generic-stable
shadowable property (resp., $C^1-$generic-stable transitive specification
property or $C^1-$generic-stable barycenter property) if and only if
$\Lambda$ is a hyperbolic basic set. In particular, $f|_\Lambda$ satisfies a $C^1-$stable
shadowable property (resp., $C^1-$stable transitive specification
property or $C^1-$stable barycenter property) if and only if
$\Lambda$ is a hyperbolic basic set.
\end{abstract}

\section{Introduction} \setlength{\parindent}{2em}

In the studies of dynamical systems, the pseudo-orbit shadowing
property usually plays an important role in the investigation of
stability theory and ergodic theory. Wen, Gan and Wen
\cite{WGW} proved that $C^1$-stably shadowable chain component is
hyperbolic. Lee, Morivasu, Sakai\cite{LMS} showed that a chain
recurrent set has $C^1-$stable shadowing property if and only if the
system satisfies both Axiom A and the no-cycle and also proved that
a chain component containing a hyperbolic periodic point $p$ has
$C^1-$stable shadowing property if and only if it is the hyperbolic
homoclinic class of $p$. Moreover, Tajbakhsh and Lee\cite{TL} proved
that a homoclinic class has $C^1-$stable shadowing property if and
only if it is hyperbolic. Recently, Sakai, Sumi and Yamamoto showed
that a closed invariant set satisfies $C^1-$stable specification
property if and only if it is a hyperbolic elementary set. Since the
specification property there naturally implies topologically mixing,
they gave a characterization of the hyperbolic (mixing) elementary
sets. Specification property is due to Bowen and Sigmund and holds
for every mixing compact set with shadowing property. In the present
paper we show that for every transitive compact set with shadowing
property, a version of transitive specification property is true.
Furthermore, we also discuss a notion called barycenter property due
to Abdenur, Bonatti, Crovisier\cite{ABC}, weaker than transitive
specification property. Here we are mainly to characterize the diffeomorphsims satisfying $C^1-$generic-stable shadowable property, transitive specification property
or barycenter property(for particular case,
$C^1-$stable shadowable property, transitive specification property
or barycenter property). More precisely, for a robustly transitive
set, it has one of above properties if and only if it is a
hyperbolic basic set.

Let $(M, d)$ denote a compact metric space and let $f :M\rightarrow
M$ be a homeomorphism. Let $\Lambda$ be a compact and $f-$invariant
set and let $f|_\Lambda$ be the restriction of $f$ on the set $\Lambda$. Now we start to introduce the notions of shadowing, specification and barycenter properties.
 A sequence
$\{y_n\}_{n=a}^b\subseteq\Lambda$ is called a $\delta$-pseudo-orbit
($\delta\geq0$) of $f$ if $d(f(y_n),y_{n+1})\leq\delta$ for every
$a\leq n\leq b$.   A system $f|_\Lambda$ is said to have the
\textit{shadowing property}  if for every $\varepsilon
>0$ there exists $\delta>0$ such that for a
given $\delta$-pseudo-orbit y=$\{y_n\}_{n=a}^b\subseteq\Lambda$, we
can find $x \in \Lambda$, which $\varepsilon$-traces y, i.e.,
$d(f^n(x), y_n)<\varepsilon$ for every $a\leq n\leq b$.
$f|_\Lambda$ is said to satisfy the \textit{transitive
specification property} if the following holds: for any
$\varepsilon>0$ there exists an increasing sequence of integers
$M_0(\varepsilon)=0<M_1(\varepsilon)< M_2(\varepsilon)<\cdots$
towards to $+\infty$ such that for any $k\geq 1,\,\, n\geq 1$, any
$k$ points $x_1,x_2,\cdots,x_k\in\Lambda$, and any integers $n_1,\,
n_2,\,\cdots,\,n_k$, there exists a point $z\in \Lambda$ and a
sequence of integers $c_1=0< c_2< \cdots< c_k$ with
$c_{j+1}-c_j-n_j\in[M_{n-1}(\varepsilon),M_n(\varepsilon)]\,(j=1,2,\cdots,k-1)$
such that $d(f^{c_{j}+i}(z), f^i(x_j))<\varepsilon,\,\,0\leq i\leq
n_j, \,\, 1\leq j\leq k.$ Now we begin to recall the barycenter
property(a little difference to \cite{ABC}).  Let $P(f|_\Lambda)$  be
the set of periodic points of $f$ in $\Lambda$. In particular,  set $P(f)=P(f|_{M})$. Given two periodic
points $p,\,q\in P(f|_\Lambda)$, we say $p,q$ have the barycenter
property, if for any $\varepsilon>0$ there exists an integer
$N=M(\varepsilon,p,q)>0$ such that for any two integers $n_1, n_2$,
there exists a point $z\in \Lambda$ and an integer $X\in[0,N]$ such
that $d(f^{i}(z), f^i(p))<\varepsilon,\,\,-n_1\leq i\leq 0, $ and
$d(f^{i+X}(z), f^i(q))<\varepsilon,\,\,0\leq i\leq n_2$.  $f|_\Lambda$ is said to
satisfy \textit{the barycenter property}  if the barycenter property holds
for any two periodic points $p,\,q\in P(f|_\Lambda)$.

Obviously the barycenter property is weaker than the transitive
specification property. The transitive specification property means
that whenever there are $k$ pieces of orbits they may be
approximated up to $\varepsilon$ by one orbit, provided that the
time for {\it switching} from the forward piece of orbit to the
afterward and the time for {\it switching back} are bounded between
two integers $M_{n-1}(\varepsilon)\leq M_{n}(\varepsilon)$, these
integers $ M_n(\varepsilon)$ being independent of the length of the
$k$ pieces of orbits.   Here this notion of \emph{transitive
specification property} is weaker than the usual \emph{(mixing)
specification property} defined by Sigmund\cite{Sig}.

Let $\Lambda$ be as before. $\Lambda$ is transitive if there is some
$x\in\Lambda$ whose forward orbit is dense in $\Lambda$. A
transitive set $\Lambda$ is trivial if it consists of a periodic
orbit.   Note that transitive specification property implies that
$\Lambda$ is topologically transitive. $\Lambda$ is locally maximal
in some neighborhood $U\subseteq M$ of $\Lambda$ if
$\Lambda= \bigcap_{k\in \mathbb{Z}} f^k(U)$. A set $\Lambda$ is a basic
set (resp. elementary set) if $\Lambda$ is locally maximal and
$f|_\Lambda$ is transitive (resp. topologically mixing). In a Baire space $X,$ we call $\mathcal{R}\subseteq X$ be a residual set, if it contains a dense $G_\delta$ set.

Let $M$ be a closed $C^{\infty}$ manifold and let $\Diff(M)$ be the
space of diffeomorphisms  of $M$ endowed with the $C^1-$topology.
Denote by $d$ the distance on $M$ induced from a Riemannian metric on
the tangent bundle $TM$. Given $f\in \Diff(M)$, denote by $\mathcal {O}(p)$ the
periodic $f-$orbit of $p\in P(f)$. If $p\in P(f)$ is a hyperbolic
saddle with period $\pi(p)>0$, then there are the local stable
manifold $W^s_\varepsilon(p)$ and the  local unstable manifold
$W^u_\varepsilon(p)$ of $p$ for some $\varepsilon=\varepsilon(p) >
0$. It is easy to see that if $d(f^n (x),f^n (p))\leq\varepsilon$
for any $n \geq 0$, then $x \in W^s_\varepsilon (p)$ (a similar
property also holds for local unstable manifold $W^u_\varepsilon
(p)$ with respect to $f^{-1}$). The stable manifold $W^s_\varepsilon
(p)$ and the unstable manifold $W^u_\varepsilon (p)$ of $p$ are
defined as usual. The dimension of the stable manifold
$W^s_\varepsilon (p)$ is called the index of $p$, and denoted by
index$(p)$.

 An $f$ invariant compact set $\Lambda$ is
robustly transitive in some neighborhood  $U$ if $\Lambda$ is locally maximal in   $U$ and there
is a $C^1$-neighborhood $\mathcal {U}(f)$ of $f$ such that for any
$g \in\mathcal {U}(f)$, $\Lambda_g(U):= \bigcap_{k\in \mathbb{Z}} g^k(U)$(called a continuation of $\Lambda_f(U)=\Lambda$) is transitive.

Now we state our main result as follows.

\begin{Thm}\label{Thm1}
Let $\Lambda$ be an $f$ invariant compact set and assume that $\Lambda$ is robustly transitive in some neighborhood  $U$. Then the following conditions are
equivalent:

(1) $f|_\Lambda$ is $C^1-$generic-stably shadowable, i.e., there is a
$C^1$-neighborhood $\mathcal {U}(f)$ of $f$ and a residual set $\mathcal{R}\subseteq \mathcal {U}(f)$ such that for any $g
\in\mathcal{R}$, $\Lambda_g(U):= \bigcap_{k\in \mathbb{Z}} g^k(U)$ has shadowing property;

(1') $f|_\Lambda$ is $C^1-$stably shadowable, i.e., there is a
$C^1$-neighborhood $\mathcal {U}(f)$ of $f$ such that for any $g
\in\mathcal {U}(f)$, $\Lambda_g(U):= \bigcap_{k\in \mathbb{Z}} g^k(U)$ has shadowing property;

(2) $f|_\Lambda$ satisfies the $C^1-$generic-stable transitive
specification property, i.e., there is a $C^1$-neighborhood
$\mathcal {U}(f)$ of $f$ and a residual set $\mathcal{R}\subseteq \mathcal {U}(f)$ such that for any $g
\in\mathcal{R}$,
$\Lambda_g(U):= \bigcap_{k\in \mathbb{Z}} g^k(U)$ has transitive specification property;

(2') $f|_\Lambda$ satisfies the $C^1-$stable transitive
specification property, i.e., there is a $C^1$-neighborhood
$\mathcal {U}(f)$ of $f$ such that for any $g \in\mathcal {U}(f)$,
$\Lambda_g(U):= \bigcap_{k\in \mathbb{Z}} g^k(U)$ has transitive specification property;

(3) $f|_\Lambda$ satisfies the $C^1-$generic-stable barycenter property,
i.e., there is a $C^1$-neighborhood $\mathcal {U}(f)$ of $f$ and a residual set $\mathcal{R}\subseteq \mathcal {U}(f)$ such that for any $g
\in\mathcal{R}$, $\Lambda_g(U):= \bigcap_{k\in \mathbb{Z}} g^k(U)$ has barycenter
property;

(3') $f|_\Lambda$ satisfies the $C^1-$stable barycenter property,
i.e., there is a $C^1$-neighborhood $\mathcal {U}(f)$ of $f$ such
that for any $g \in\mathcal {U}(f)$, $\Lambda_g(U):= \bigcap_{k\in \mathbb{Z}} g^k(U)$ has barycenter
property;

(4) there is a $C^1$-neighborhood $\mathcal {U}(f)$ of $f$ and a residual set $\mathcal{R}\subseteq \mathcal {U}(f)$ such that for any $g
\in\mathcal{R}$, any two periodic hyperbolic saddles
$p,\,q\in\Lambda_g(U):= \bigcap_{k\in \mathbb{Z}} g^k(U),$ $W^s(\mathcal {O}(p))\cap W^u(\mathcal
{O}(q))\neq\emptyset;$

(4') there is a $C^1$-neighborhood $\mathcal {U}(f)$ of $f$ such that
for any $g \in\mathcal {U}(f)$, any two periodic hyperbolic saddles
$p,\,q\in\Lambda_g(U):= \bigcap_{k\in \mathbb{Z}} g^k(U),$ $W^s(\mathcal {O}(p))\cap W^u(\mathcal
{O}(q))\neq\emptyset;$

(5) there is a $C^1$-neighborhood $\mathcal {U}(f)$ of $f$ such that
for any $g \in\mathcal {U}(f)$, any  periodic point of
$\Lambda_g(U):= \bigcap_{k\in \mathbb{Z}} g^k(U)$ is hyperbolic and has the same index;

(6) there is a $C^1$-neighborhood $\mathcal {U}(f)$ of $f$ such that
for any $g \in\mathcal {U}(f)$, $\Lambda_g(U):= \bigcap_{k\in \mathbb{Z}} g^k(U)$ is a hyperbolic basic
set.
\end{Thm}

{\bf Remarks.} 1. The result in \cite{SSY} is a particular case of $(2')$ in Theorem \ref{Thm1} (being as the mixing case), since the specification property assumed on $\Lambda_g(U)$ in \cite{SSY} naturally implies that $\Lambda_g(U)$ is topologically mixing(i.e., $\Lambda$ is robustly mixing). And we point out that the codition $(2)$ in Theorem \ref{Thm1} can be also as a generalization of the result\cite{SSY}, since $(2)$ is weaker than $(2')$.

2. From \cite{SSY} the
set of transitive Anosov diffeomorphisms is a characterization of the set of
diffeomorphisms satisfying $C^1-$stable mixing specification
property. Here by Theorem \ref{Thm1} it is also a
characterization of the set of diffeomorphisms satisfying $C^1-$stable (or generic-stable) transitive specification
property. Furthermore, if $M$ is topologically transitive, the set
of Anosov diffeomorphisms is also a characterization of
$C^1-$stable (or generic-stable) shadowable property or $C^1-$stable (or generic-stable) barycenter property.

\medskip
 The equivalence of (5) and (6) is due to Ma\~{n}$\acute{\text{e}}$ \cite{Ma} and by
using this Sakai, Sumi and Yamamoto \cite{SSY} proved that
$f|_{\Lambda_f(U)}$ satisfies the $C^1-$stable specification
property (mixing case) if and only if $\Lambda$ is a hyperbolic
elementary set. Actually, it is essentially proved
$(4')\Rightarrow(5)$ in \cite{SSY}. More precisely, $(4')$ implies any
two hyperbolic periodic saddles $p,\,q\in\Lambda_g(U)$ have the same
index (see the proof of Lemma
 2.2 in \cite{SSY}) and the later implies all periodic points in $\Lambda_g(U)$ are
 hyperbolic (see Lemma 2.4 in \cite{SSY}). $(6)\Rightarrow(1')$,
 $(a')\Rightarrow(a)(a=1,2,3,4)$ and $(2)\Rightarrow(3)$
 are obvious and thus it is enough to show
 $(1)\Rightarrow(2)$, $(3)\Rightarrow(4)$ and $(4)\Rightarrow(5)$.

\section {\bf Proof of our main theorem}

To prove $(1)\Rightarrow(2)$, $(3)\Rightarrow(4)$ and $(4)\Rightarrow(5)$ in our main
theorem, we divide into three lemmas. Firstly, we show a general lemma
which implies $(1)\Rightarrow(2)$.

\begin{Lem}\label{Lem5} Let $f:M\rightarrow M$ be a homeomorphism on a
compact metric space $M$ and let $\Lambda$ be a transitive
$f-$invariant subset. If $f|_\Lambda$  satisfies the shadowing
property, then $\Lambda$ satisfies the transitive specification
property.
\end{Lem}
Remark: If we assume $\Lambda$ is mixing in Lemma \ref{Lem5}, then
$(f,\Lambda)$ satisfies the mixing
specification property of Sigmund\cite{Sig}.\\

{\bf Proof of Lemma \ref{Lem5}}  For any $\varepsilon>0$, by
shadowing property there exists $\delta>0$ such that any
$\delta-$pseudo orbit in $\Lambda$ can be $\varepsilon$ shadowed by
a true orbit in $\Lambda$.

Take and fix for $\Lambda$ a finite cover
$\alpha=\{U_1,U_2,\cdots,U_{r_0}\}$ by nonempty open balls $U_i$ in
$\Lambda$ satisfying $diam(U_i) <\delta$, $i=1,\,2,\,\cdots,\,r_0$.
Since $\Lambda$ is transitive, for any $i,\,j=1,\,2,\,\cdots,\,r_0,$
there exist a positive integer $X^{(1)}_{i,\,j}$ such that
$$f^{-X^{(1)}_{i,\,j}}(U_i)\cap
  U_j\neq \emptyset.$$ Let
$$M_{1}=max_{1\leq
i\neq j\leq r_0}X^{(1)}_{i,\,j}.$$ Similarly for any
$i,\,j=1,\,2,\,\cdots,\,r_0,$ we can take a positive integer
$X^{(2)}_{i,\,j}\geq M_1$ such that $$f^{-X^{(2)}_{i,\,j}}(U_i)\cap
  U_j\neq \emptyset.$$ Let
$$M_{2}=max_{1\leq
i\neq j\leq r_0}X^{(2)}_{i,\,j}.$$ By induction for any
$i,\,j=1,\,2,\,\cdots,\,r_0,$ there is a sequence of increasing
integers $1\leq X^{(1)}_{i,\,j}< X^{(2)}_{i,\,j}<\cdots<
X^{(n)}_{i,\,j}<\cdots<+\infty$ such that
$$f^{-X^{(n)}_{i,\,j}}(U_i)\cap
  U_j\neq \emptyset$$ and $$ X^{(n)}_{i,\,j}\geq M_{n-1},$$ where
$$M_{n-1}=max_{1\leq
i\neq j\leq r_0}X^{(n-1)}_{i,\,j}.$$ Setting $M_0=0,$ clearly
$\{M_n\}_{n\geq0}$ is an increasing sequence towards to $+\infty.$

Now let us consider a given  sequence of points
$x_1,\,x_2,\,\cdots,\,x_k\in\,\Lambda,$ and a sequence of positive
numbers ${n_1,\,n_2,\,\cdots,\,n_k}$. Take and fix
$U_{i_0},\,U_{i_1}\in\alpha$ so that $x_i\in
U_{i_0},\,\,f^{n_i}(x_i)\in U_{i_1},i=1,\,2,\,\cdots,\,k.$ Fixing an
integer $n\geq1,$ take $y_i\in U_{i_1}$ such that
$f^{X^{(n)}_{(i+1)_0,\,i_1}}(y_i)\in U_{(i+1)_0}$ for
$i=1,\,2,\,\cdots,\,k-1.$ Take $y_k\in U_{k_1}$ such that
$f^{X^{(n)}_{(k+1)_0,\,k_1}}(y_k)\in U_{1_0}$. Thus we get a periodic
 $\delta$-pseudoorbit in $\Lambda$:
$$\{f^t(x_1)\}_{t=0}^{n_1}\cup
\{f^t(y_1)\}_{t=0}^{X^{(n)}_{2_0,1_1}}\cup \{f^t(x_2)\}_{t=0}^{n_2}
\cup\cdots\cup\{f^t(x_k)\}_{t=0}^{n_k}
\cup\{f^t(y_k)\}_{t=0}^{X^{(n)}_{(k+1)_0,\,k_1}}.$$  Hence there
exists a point $z\in \Lambda$ $\varepsilon$-shadowing the above
sequence. More precisely,
$$d(f^{c_{i-1}+j}(z),f^j(x_i))<\varepsilon,\,\,j=0,\,1,\,\cdots,\,n_i,\,\,\,i=1,\,2,\,\cdots,\,k,$$
where $c_i$ is defined
 as follows:
 $$c_i=\begin{cases}
 \,\,0,&\text{for }\,\,i=0\\
 \,\,\sum_{j=1}^{i}[n_j+X^{(n)}_{(j+1)_0,\,j_1}],\,\,&\text{for}\,\,\,\,i=1,\,2,\,\cdots,\,k.\\
\end{cases}
 $$ \hfill $\Box$
\bigskip

Secondly, we prove a lemma about the relationship of homoclinic
related property and barycenter property, which deduces
$(3)\Leftrightarrow(4)$ of our main theorem.

\begin{Lem}\label{Lem1} Let $f:M\rightarrow M$ be a diffeomorphism on a
compact manifold $M$. Then for two hyperbolic periodic points
$p,\,q\in P(f),$  $p,\,q$ have the barycenter property
$\Leftrightarrow$ $W^u(\mathcal {O}(p))\cap W^s(\mathcal
{O}(q))\neq\emptyset$.
\end{Lem}

 {\bf Proof of Lemma \ref{Lem1}}

$"\Rightarrow":$ Let $\varepsilon(p)$ and $\varepsilon(q) > 0$ be as
before with respect to $p$ and $q$. Take
$\varepsilon=min\{\varepsilon(p),\,\varepsilon(q)\}$, and let $N =
N(\varepsilon, p, q)
> 0$ be the number of barycenter property. For any $n\geq0$, by barycenter property there is $z_n\in \Lambda$
and an integer $X_n\in[0,N] $ such that \\$(i).\,\,\,\,\,\, d(f^j
(z_n ), f^j (p))\leq\varepsilon\,\,\text{ for} \,\,-n\leq j \leq0,$\\
$(ii).\,\,\,\,\,\, d(f^{j+X_n} (z_n), f^j
(q))\leq\varepsilon\,\,\text{ for} \,\,0\leq j \leq n.$

Take a subsequence $\{n_k\}$ such that $n_k\rightarrow\infty$ and
$X_{n_k}\equiv X$ for some fixed integer $X\in[0,N].$ Let $z
=\lim_{k\rightarrow+\infty} z_{n_k}$ by taking a subsequence again
if necessary. By (i) and (ii) one has $z\in
W^u_{\varepsilon(p)}(p)\subseteq W^u(p)$ and $f^X(z)\in
W^s_{\varepsilon(q)}(q)\subseteq W^u(q)$.\\

$"\Leftarrow":$ If $W^u(\mathcal {O}(p))\cap W^s(\mathcal
{O}(q))\neq\emptyset$, then we can take $z\in W^u(\mathcal
{O}(p))\cap W^s(\mathcal {O}(q)).$ So for any $\varepsilon>0,$ there
is $N_1=N_1(\varepsilon,p,q)>0$ such that $d(f^j(z),f^j(p))<\varepsilon$
for all $j\leq-N_1$ and $d(f^j(z),f^j(q))<\varepsilon$ for all $j\geq
N_1$. Moreover, we can assume $N_1$ to be  a common multiple of the
period of $p$ and $q.$ Put $x=f^{-N_1}(z)$ and let $N=2N_1.$ Then
for any two integers $n_1,\,n_2,$ $x$ is needed for barycenter
property, i.e., $d(f^j(x),f^j(p))=d(f^{j-N_1}(z),f^{j-N_1}(p))<\varepsilon$
for all $-n_1\leq j\leq0$ and
$d(f^{j+N}(x),f^j(q))=d(f^{j+N_1}(z),f^{j+N_1}(q))<\varepsilon$ for all
$0\leq j\leq n_2$. \hfill $\Box$
\bigskip

Before proving $(4)\Rightarrow(5)$, we state a lemma which says that $(4)\Rightarrow$ all hyperbolic points have the same index. A diffeomorphism $f$ is said to be $Kupka-Smale$ if the periodic points of $f$ are hyperbolic and for any two periodic points $p,q$ of $f$, $W^s(p)$ is transversal to $W^u(q).$ It is well known that the set of Kupka-Samle diffeomorphisms is $C^1$-residual in $\Diff(M)$ (see\cite{Ro}). Note that if $\mathcal{U}$ is an open set of  $\Diff(M)$, then the set of Kupka-Samle diffeomorphisms restricted in $\mathcal{U}$ is still $C^1$-residual in $\mathcal{U}$.

\begin{Lem}\label{Lem2} Let $f:M\rightarrow M$ be a diffeomorphism on a
compact manifold $M$. Then condition $(4)$ in Theorem \ref{Thm1} implies that for any two hyperbolic saddles
$p,\,q\in \Lambda_g(U)\cap P(g)$ with respect to $g\in \mathcal{U}(f),$   $index(p)=index(q)$.
\end{Lem}

{\bf Proof} This proof is an adaption of Lemma 2.2 in \cite{SSY}. Let $\mathcal{U}(f)$ be as in condition $(4)$ of Theorem \ref{Thm1}. Fix a $g\in \mathcal{U}(f),$ and let $p,q\in\Lambda_g(U)\cap P(g)$ be hyperbolic saddles. Then
there is a $C^1-$neighborhood $\mathcal{V}(g)\subseteq \mathcal{U}(f)$ such that for any $\varphi\in  \mathcal{V}(g)$,
there is continuations $p_\varphi$ and $q_\varphi$ (of $p$ and $q$) in $\Lambda_\varphi(U)$, respectively (Since $\Lambda_\varphi(U)=\Lambda\subset int U$, we can assume that $\Lambda_g(U)\subset int U$ for any $g\in \mathcal{U}(f)$ reducing $ \mathcal{U}(f)$ if necessary).

By contradiction, if $index(p)<index(q)$(the other case is similar), then we have $$dim W^s(p,g) + dim W^u(q,g) < dim M,$$ where $W^s(p,g)$ and $W^u(q,g)$ are the stable and unstable manifold of $p$ and $q$ with respect to $g$. Since the intersection of two residual sets is still residual, then the set of diffeomorphisms restricted in $\mathcal{V}(g)$ satisfying not only Kupka-Smale but also condition (4) of Theorem \ref{Thm1} is still residual in $\mathcal{V}(g)$. Take such a diffeomorphism $\varphi\in\mathcal{V}(g)$. Then $$W^s(p_\varphi,\varphi) \cap W^u(q_\varphi,\varphi)=\varnothing,$$
since $dim W^s(p,g)=dim W^s(p_\varphi,\varphi) $ and $dim W^u(q,g)=dim W^u(q_\varphi,\varphi).$ On the other hand, since $\varphi$ is a diffeomorphism satisfying condition (4) of Theorem \ref{Thm1}, then $$W^s(p_\varphi,\varphi) \cap W^u(q_\varphi,\varphi)\neq\varnothing.$$ This is a contradiction.
 \hfill $\Box$

\bigskip
{\bf End of proof of (4)$\Rightarrow$(5):} By lemma \ref{Lem2} we only need to prove that every periodic point
$p\in \Lambda_g(U)$ of $g\in \mathcal{U}(f)$ is hyperbolic. By contradiction, suppose that $p\in \Lambda_g(U)$ of $g\in \mathcal{U}(f)$ is not hyperbolic. Then by Lemma 2.4 in \cite{SSY}, there is $\varphi\in \mathcal{U}(f)$ possessing hyperbolic points $q_1$ and $q_2$ in $\Lambda_{\varphi}(U)$ with different indices. This is a contradiction to Lemma \ref{Lem2}. \hfill $\Box$

\section{One remark for volume-preserving version.}

Let  $\omega$ be a volume measure on the smooth compact manifold $M$ and $\Diff_\omega(M)$ be the space of diffeomorphisms preserving $\omega$. We point out the statements in Theorem \ref{Thm1} can be changed for the volume-preserving diffeomorphisms, since all the main techniques can be replaced by the ones of
volume-preserving version. I.e., the
set of transitive volume-preserving Anosov diffeomorphisms is a characterization of the set of volume-preserving
diffeomorphisms satisfying 
$C^1-$stable (or generic-stable) shadowable property, $C^1-$stable (or generic-stable) transitive or mixing specification
property or $C^1-$stable (or generic-stable) barycenter property.  Let's explain it more precisely as follows. The equivalence of  condition (5) and (6) (see \cite{Ma}) in Theorem \ref{Thm1} can be replaced by the recent result in \cite{AC}, Frank's lemma\cite{Fra} (important to prove Lemma 2.4 in \cite{SSY} which is needed in our proof of (4)$\Rightarrow$(5))  can be replaced by the pasting lemma for volume-preserving systems(see \cite{AM}) and Kupka-Smale property for volume-preserving case can be found in \cite{Ro1}. In particular, we  point out that we need not assume the robust transitivity of $M$ when we prove the result that one volume-preserving diffeomorphism satisfying
$C^1-$stable (or generic-stable) shadowable property  is Anosov,  since generic   volume-preserving diffeomorphisms are transitive from \cite{BC}(and so that by Lemma \ref{Lem5} generic-stable shadowing implies generic-stable transitive specification for volume-preserving). Moreover, we note that volume-preserving Anosov diffeomorphisms are always transitive from the viewpoint of structurally stable property of Anosov systems and the transitivity of generic   volume-preserving diffeomorphisms\cite{BC}, since every volume-preserving Anosov diffeomorphism has a topologically conjugated diffeomorphism arbitrarily nearby which can also be chosen transitive from \cite{BC} and transitivity is an invariant property under conjugation. So the statements above for volume-preserving case can be  directly as a characterization of (not necessarily adding $``$transitive") volume-preserving Anosov diffeomorphisms.
\bigskip

{\bf Acknowledgement} It is a great pleasure to thank the referee for valuable comments and suggestions.

\end{document}